\documentclass[letterpaper, 10 pt, conference]{ieeeconf}  

\IEEEoverridecommandlockouts                              
\usepackage{amsmath, amssymb, graphicx, bm}
\usepackage{algorithm}
\usepackage{algpseudocode}
\usepackage{subcaption, float}
\usepackage{hyperref}
\usepackage{cleveref}
\usepackage{placeins}
\usepackage[absolute]{textpos}

\usepackage{xcolor}

\newcommand{\nor}[1]{\left\| #1 \right\|} 
\newcommand{\LRp}[1]{\left( #1 \right)} 

\newcommand{\mc}[1]{\mathcal{#1}} 
\newcommand{\mbb}[1]{\mathbb{#1}} 

\title{\LARGE \bf
Rendezvous Planning from Sparse Observations of Optimally Controlled Targets
}

\author{Thomas A. Scott, Lukas Taus, Yen-Hsi Richard Tsai, Tan Bui-Thanh, Justin G.R. Delva
\thanks{This work was supported in part by Lockheed Martin.}
\thanks{Thomas A. Scott is with the Dept of Aerospace Engineering and Engineering Mechanics, UT Austin, USA (\tt\small thomas.scott@utexas.edu)}
\thanks{Lukas Taus is with the Oden Institute, UT Austin, USA (\tt\small l.taus@utexas.edu)}
\thanks{Yen-Hsi Richard Tsai is with Dept of Mathematics and Oden Institute, UT Austin, USA (\tt\small ytsai@math.utexas.edu)}
\thanks{Tan Bui-Thanh is with the Dept of Aerospace Engineering and Engineering Mechanics and the Oden Institute, UT Austin (\tt\small tanbui@utexas.edu)}
\thanks{Justin G.R. Delva is with Lockheed Martin, USA (\tt\small justin.delva@lmco.com)}
}

\begin{document}

\maketitle
\thispagestyle{empty}
\pagestyle{empty}

\setlength{\TPHorizModule}{\textwidth}
\setlength{\TPVertModule}{\textwidth}
\begin{textblock}{0.7}(0.25,1.49)
\footnotesize \centering This work has been submitted to the IEEE for possible publication. Copyright may be transferred without notice, after which this version may no longer be accessible.
\end{textblock}

\begin{abstract}

We develop a probabilistic framework for \emph{rendezvous planning}: given sparse, noisy observations of a fast-moving target, plan rendezvous spatiotemporal coordinates for a set of significantly slower seeking agents. The unknown target trajectory is estimated under uncertain dynamics using a filtering approach that combines a kernel-based maximum a posteriori estimation with Gaussian process correction, producing a mixture over trajectory hypotheses.
This estimate is used to select spatiotemporal rendezvous points that maximize the probability of successful rendezvous. Points are chosen sequentially by greedily minimizing failure probability in the current belief space, which is updated after each step by conditioning on unsuccessful rendezvous attempts. We show that the failure-conditioned update correctly captures the posterior belief for subsequent decisions, ensuring that each step in the greedy sequence is informed by a statistically consistent representation of the remaining search space, and derive the corresponding Bayesian updates incorporating temporal correlations intrinsic to the trajectory model.
This result provides a systematic framework for planning under uncertainty in applications of autonomous rendezvous such as unmanned aerial vehicle refueling, spacecraft servicing, autonomous surface vessel operations, search and rescue missions, and missile defense.
In each, the motion of the target entity can be modeled using a system of differential equations undergoing optimal control for a chosen objective, in our example case Hamilton--Jacobi--Bellman solutions for minimum arrival time of a Dubins car with uncertain turning radius and destination.

\end{abstract}


\section{INTRODUCTION}
The rapid advancement of autonomous systems has shifted the paradigm of aerial and maritime monitoring from reactive tracking to the analysis of goal-oriented behavior. As agents increasingly operate with high degrees of onboard intelligence, predicting their future trajectories requires models that account for the underlying optimality of their decision-making processes. In many real-world scenarios, however, this predictive task is characterized by an ``information gap'' caused by intermittent observations, uncertainty about the target agent's exact kinematics, and environmental constraints. This challenge is further compounded when the target agent possesses a significant velocity advantage over the responding pursuer. 
Examples include an unmanned aerial vehicle (UAV) attempting to deliver a payload to a moving recipient such as for mid-air refueling~\cite{moreira2019, pliska2024}, spacecraft servicing~\cite{mohammad2024}, space debris removal~\cite{maestrini2022}, search and rescue missions involving UAVs~\cite{doschl2025}, autonomous surface vessel operations~\cite{mirabito2017}, and a pursuer attempting to rendezvous with an evader such as in a defensive missile scenario~\cite{moreira2019}. In such settings, this problem is often addressed through Kalman filtering methods and proportional guidance or model predictive control~\cite{moreira2019, pliska2024}. However, while Kalman filtering and MPC provide a foundational approach, they do not inherently leverage the reachability properties of the target, which are essential for guaranteeing rendezvous when observations are both rare and noisy.

The current state estimate is propagated under an assumed governing equation with inferred parameters to predict future target states — typically via interacting Kalman filters with tailored linear dynamics \cite{pliska2024, genovese2001}. However, this only addresses the unknown parameters and state but not potential inaccuracy of the dynamics. The exact trajectory is uncertain due to noise in perception, environmental variability, or behavioral deviation from the assumed dynamics undergoing optimal control. We model this uncertainty through a Gaussian process (GP), i.e. a stochastic process with a parameterized mean and covariance. After estimating a maximum a posteriori (MAP) trajectory, we propose using Gaussian process regression (GPR) with the MAP trajectory as the prior to provide the stochastic model. This can be viewed as a filtering step followed by a data assimilation step. We repeat this across potential parameter variations and through Bayesian inference arrive at a model of the uncertain trajectory.
This is conceptually similar to existing data assimilation methods that utilize GPs such as the Kennedy-O’Hagan method \cite{kennedy_ohagan2002} and the subsequent improvements \cite{arcones2025}. However, these existing methods are costly, making them only suitable for offline data assimilation and not for real-time applications such as tracking for rendezvous.

Formally, we consider the problem of planning spatiotemporal rendezvous points that maximize the probability of successful rendezvous — a setting we term \emph{rendezvous planning}.
Traditionally, proportional guidance is used to follow the target entity, but this fails when the pursuer is significantly slower than the target: reactive tail-chasing becomes geometrically infeasible, making predictive planning essential rather than merely beneficial. Planning for rendezvous overcomes this limitation while accommodating additional operational constraints such as a desired contact angle.

To guarantee that the computed rendezvous points are attainable in the required time, we utilize the reachable region $\mathcal{R}$. This set defines the spatiotemporal coordinates attainable by the pursuer given its specific dynamics and initial configuration. By incorporating $\mathcal{R}$ directly into our algorithm, we ensure that proposed rendezvous points not only lead to high likelihood of contact but are also physically feasible. This framework allows us to enforce complex operational constraints, such as specified contact angles (e.g., side-on approaches) and restricted initial configurations (e.g., fixed launch headings), which are difficult to incorporate into traditional reactive guidance laws.

We demonstrate that our method is more effective than Kalman filter paired with proportional guidance in the sparse observation regime, a common state of the art strategy for non-cooperative rendezvous problems. In this regime, a linear dynamics model (prior prediction of Kalman filter step) is very inaccurate at predicting the true, nonlinear dynamics and the sparse quantity of data is insufficient to correct for this. In contrast, MAP estimate models the nonlinear dynamics directly and the subsequent GPR step can correct for model inadequacy and provide uncertainty quantification (UQ). In the examples considered here, the proposed estimation pipeline yields target-state predictions sufficiently accurate to support rendezvous planning from sparse, highly noisy observations.

\section{PROBLEM FORMULATION}

The target entity is assumed to follow a trajectory that is approximately optimal with respect to a minimum travel-time objective toward an unknown goal, and can therefore be modeled by an ordinary differential equation (ODE) with unknown parameters. Since both the destination and the motion parameters of the target entity are unknown, multiple approximately optimal candidate trajectories must be computed. An efficient way to generate these trajectories is through a Hamilton–Jacobi formulation of optimal control, which enables the extraction of multiple trajectories to a destination from a single solve.
The state, i.e. the trajectory function evaluated at time $t\in\mbb{R}$, is modeled as the state $z(t)\in\mc{Z}=\mbb{R}^{d_z}$ with ODE constraint $\partial_t z(t) = \Phi(t,z(t);p)$, where $p\in\mc{P}=\mc{P}_1\times...\times\mc{P}_{d_p}$ denotes the parameters.
The estimated state of the target is periodically refined through noisy observations $y\in \mc{Y}$, that at least includes the position $x$ but may also include other variables dependent on the state $z$.
The observations are assumed to be i.i.d. normal after a sensing function i.e. $y_i\sim \mc{N}(\text{sensor}_i(z(t)), \sigma_{s,i}^2)$ where $\text{sensor}_i:\mc{Z} \rightarrow \mbb{R}$.
All observations $y$ are used in the MAP estimation step, while only position observations are used in the GPR step.

\subsection{Dubins Car}

We consider the Dubins car model to test our methodology. This model, originally developed for ground vehicles with a constrained turning radius, has also been widely adopted for UAVs \cite{burns2007} and autonomous surface vessels \cite{mirabito2017}, as ships and boats are similarly subject to minimum turning radius constraints. The model represents a target entity that follows a minimum arrival time trajectory to some candidate destination region $\Gamma_i$. We consider the case where a set of possible candidate destinations $\Gamma_1 ... \Gamma_{N_\Gamma}$ are known. The Dubins car dynamics represent a vehicle following a two-dimensional trajectory with position $x=(x_1,x_2)$ and orientation $\theta$ as the state $z=(x_1,x_2,\theta)$ maintaining a constant speed $v \in \mbb{R}^+$ and whose path is limited by a minimum turning radius $\rho \in \mbb{R}^+$.
We directly fully observe the state in this example i.e. $\mc{Y}=\mc{Z}$ and $\text{sensor}_i(z)=z_i$.
We treat $\rho$ as an unknown parameter along with the index $i$ corresponding to $\Gamma_i$, i.e. $p=(\rho,i)$. Thus $\rho$ captures uncertainty in the ODE dynamics, while $i$ captures uncertainty in the terminal condition of the HJB equation — both are inferred from observations.
Formally target vehicle dynamics and control constraint are
\begin{equation}
\dot{x}_1 = v\cos \theta, \hspace{4pt} \dot{x}_2 = v\sin \theta, \hspace{4pt} \dot{\theta} = \frac{\alpha v}{\rho}, \quad \alpha\in[-1,1],
\label{eq:car_ode}
\end{equation}
where $\alpha$ is the control parameter whose space enforces the turning radius constraint.
The unknown parameter considered is thus $p = (\rho,i) \in \mbb{R}^+ \times (1 \dots N_\Gamma)$.
The optimal behavior, a path to a given destination $\Gamma_i$ in minimum time, can be modeled by a Hamilton–Jacobi–Bellman (HJB) equation \cite{const_curve} that describes time-optimal or cost-minimizing motion under known environmental constraints:
\begin{equation} \label{eq:HJB_PDE}
0 = \cos(\theta) u_{x_1} + \sin(\theta) u_{x_2} - \frac{1}{\rho}\vert u_\theta \vert + \frac{1}{v}, \quad u \equiv 0 \text{ in } \Gamma_i.
\end{equation}
Plugging in the optimal control to the dynamics we get the final parameterized ODE constraint for Dubins car,
\begin{equation}
    \Phi(x_1,x_2,\theta;\rho,i) = \begin{bmatrix}
    v \cos \theta \\
    v \sin \theta \\
    -\frac{v}{\rho} \text{sgn}(\frac{\partial u}{\partial \theta} (x_1,x_2,\theta;\rho,\Gamma_i))
    \end{bmatrix}
    \label{eq:car_ode_constraint}
\end{equation}

To solve Eq.~\ref{eq:HJB_PDE}, we employ an upwind finite difference scheme as described in \cite{const_curve}. However, for more complex dynamical systems, particularly those with high-dimensional state spaces, this approach becomes intractable. In such cases, a machine learning method for efficiently solving high-dimensional Hamilton–Jacobi equations has been proposed in \cite{NN_LxF}.

\subsection{Proportional Guidance Baseline}
For comparison with existing methods, we evaluate our proposed framework against a Proportional Guidance (PG) law coupled with a Kalman Filter. PG is fundamentally reactive and assumes a favorable velocity ratio. In the regime where the target moves significantly faster than the pursuer ($v_e \gg v_p$), the pursuer cannot rely on closing velocity to compensate for estimation lags or sparse data. Furthermore, while more sophisticated filtering techniques such as Interacting Multiple Model (IMM) filters could improve state estimation, they do not inherently address the reachability constraints imposed by the target's optimal control policy. By using PG as a baseline, we demonstrate that traditional reactive guidance is insufficient for rendezvous under these constraints, necessitating the predictive, reachability-aware approach developed in this work.
The PG law can be described as
\begin{align}
& \dot{x}_1 = v\cos \theta, \quad \dot{x}_2 = v\sin \theta, \label{eq:prop_ode} \\
& \dot{\theta} = C_p \text{atan2}(\Delta x_2 \cos \theta - \Delta x_1 \sin \theta, \Delta x_1 \cos \theta + \Delta x_2 \sin \theta) \nonumber,
\end{align}
where $\text{atan2}(y,x)$ is the two-argument arctangent, $(\Delta x_1, \Delta x_2)$ is the difference between the Kalman filter state estimate of the target and the current pursuer state, and the gain $C_p=\frac{v}{\pi \rho}$ links the proportional constant to the minimum turning radius.
Unlike the Dubins car, which immediately takes the sharpest available turn to achieve the desired heading, proportional guidance does so only in the degenerate case where the pursuer is heading directly away from the target.
The Kalman filter uses a constant-acceleration model, treating jerk as Gaussian noise.

\section{STOCHASTIC TRAJECTORY ESTIMATION}

\subsection{ODE-Constrained MAP Estimation}
Our method begins with an optimal recovery step to compute a MAP trajectory of the target consistent with the assumed ODE dynamics.
Specifically, we seek the minimum-norm solution subject to the ODE constraint, with a squared data loss over $n_y$ observations at times $\bm{t}^y$:
\begin{align}
    \underset{z}{\min} \hspace{4pt} &\sum_{i=1}^{d_z} \nor{z_i}_\mc{H} + \sum_{i=1}^{d_y} \frac{1}{\beta_i^2} \sum_{j=1}^{n_y} (\bm{y}_{i,j} - \text{sensor}_i(z(t_j)))^2 \nonumber \\
    & s.t. \hspace{4pt} \partial_t z = \Phi(\cdot,z;p)
    \quad \forall t \in [0,T)
    \label{eq:gp1_odes}
\end{align}
where we use a reproducing kernel Hilbert space (RKHS) $\mc{H}$ defined by a chosen kernel $k:\mbb{R}^2\rightarrow\mbb{R}$
Following the general method developed by \cite{chen_solving_2021}, we can reformulate \cref{eq:gp1_odes} as a finite-dimensional optimization problem via the kernel trick.
After discretizing — yielding a problem equivalent to standard GPR when constraints are removed and $\text{sensor}_i$ is the identity —
\begin{align}
    \underset{\bm{w}}{\min} \quad &
    \sum_{i=1}^{d_z}
    \bm{w}_i^\top \Theta_i^{-1} \bm{w}_i
    + \sum_{i=1}^{d_y} \frac{1}{\beta_i^2} \sum_{j=1}^{n_y} (\bm{y}_{i,j} - \text{sensor}_i(\bm{w}_{y,j}))^2
    \nonumber
    \\
    & s.t. \quad \partial_t z(t) = \Phi(t,z(t);p) \quad \forall t \in \bm{t}^\Phi
    \label{eq:gp1_finite}
\end{align}
where $\bm{w}$ contains pointwise values of the state $z$ and its time derivatives to be regressed by the kernel. Subscripts $0$ and $1$ denote function values and first derivatives at the ODE collocation points $\bm{t}^{\Phi}$; the $y$ subscript denotes values at the $n_y$ observation times. The matrices $\Theta_i$ are the corresponding Gram matrices.
Explicitly, $\Theta_i$ corresponding to $\bm{w}_i=[\bm{w}_{i,y}^\top,\bm{w}_{i,0}^\top, \bm{w}_{i,1}^\top]^\top$ is the block matrix
\begin{equation}
    \Theta_i =
    \begin{bmatrix}
        k(\bm{t},\bm{t}) & (\partial'_t k)(\bm{t},\bm{t}^\Phi) \\
        (\partial_t k)(\bm{t}^\Phi,\bm{t}) & (\partial_t \partial'_t k)(\bm{t}^\Phi,\bm{t}^\Phi) \\
    \end{bmatrix},
\end{equation}
where $\partial'$ denotes the derivative w.r.t. the second argument, $k(\bm{t},\bm{t})$ is the matrix $K_{i,j}=k(t_i,t_j)$, and $\bm{t} = [(\bm{t}^y)^\top, (\bm{t}^\Phi)^\top]^\top$.
Meanwhile, the diagonal perturbation matrix that we will introduce for numerical stability is
\begin{equation}
    D_i^\eta =
    \begin{bmatrix}
        \eta I_{\# \bm{t}} & \\
        & \eta \frac{k_i(0,0)}{(\partial_t \partial'_t k_i)(0,0)} I_{\# \bm{t}^\Phi} \\
    \end{bmatrix},
\end{equation}
where $0 \leq \eta \ll 1$ is the typical nugget parameter introduced to make the matrix well conditioned enough to factor numerically.
After forcing constraint satisfaction in the construction of $\bm{w}$ by simply setting the value of $\bm{w}_{i,1}$ to the RHS of the $\Phi$ in Equation~\eqref{eq:gp1_finite},
\begin{equation}
\begin{aligned}
    & \min_{\substack{\bm{w}_{1,0},\bm{w}_{1,y}, \dots, \\ \bm{w}_{d_z,0},\bm{w}_{d_z,y}}} \quad 
    \sum_{i=1}^{d_z} \left( \bm{w}_i^\top (\Theta_i + D_i^\eta)^{-1} \bm{w}_i \right) \\
    & \quad + \sum_{i=1}^{d_y} \left( \frac{1}{\beta_i^2} \sum_{j=1}^{n_y} (\bm{y}_{i,j} - \text{sensor}_i(\bm{w}_{y,j}))^2 \right).
    \label{eq:gp1_final}
\end{aligned}
\end{equation}

After optimization, the MAP solution can be evaluated by computing the coefficients $\bm{\alpha}_i := (\Theta_i + D_i^\eta)^{-1} \bm{w}_i$ then
\begin{equation}
    z_i^*(t) = \sum_{j=1}^{\# \bm{t}} \alpha_{i,j} k(t,\bm{t}_j) + \sum_{j=\# \bm{t}+1}^{\# \bm{t}^{\Phi}} \alpha_{i,j} (\partial_t k)(t,\bm{t}^{\Phi}_j).
    \label{eq:map_eval}
\end{equation}

\subsection{Gaussian Process Trajectory Model}
We model both the trajectory $x(t)$ and parameters $p$ stochastically, denoting all probability density functions as $f$.
The target trajectory is modeled as a stochastic process $X_t$; conditioned on parameters $p$, we write $X_t|p$ for the corresponding process with conditional density $f(x|t,p)$.
Each component of state $z$ in each $X_t|p$ is modeled as an a priori independent GP, with the MAP position estimate $x^*:=\text{pos}(z^*)(t)$ as the prior mean and the kernels $k_i$ from the optimal recovery step as the prior covariance:
\begin{equation}
    X_t|p \sim \mathcal{GP}\LRp{
        x^*,
        \begin{bmatrix}
            k_1^\sharp(t,t') & & \\
            & \ddots & \\
            & & k_{d_x}^\sharp(t,t')
        \end{bmatrix}
    }.
\end{equation}
Furthermore, we assume each component of the parameter is a priori independent.
We will use a Monte Carlo approximation for marginalization by sampling a finite set $\mc{P}_\text{MC} \subset \mc{P}$.

\subsubsection{GPR for Uncertainty Quantification and Bias Correction}
The second step is standard GPR \cite{rasmussen:GaussianProcesses} with the MAP estimate as the prior mean, serving dual purposes: uncertainty quantification and correction of residual model-data discrepancies.
Starting from the marginal likelihood of $X_t|p$ conditioned on $n_y$ position observations $\bm{y}_x$ at times $\bm{t}^y$ given parameter $p$, we apply Bayes' theorem and marginalize over $p$ to obtain the posterior density of $x$.
For brevity, the MAP model and remaining GP hyperparameters are suppressed in the conditional notation but remain implicit.
In log form, the density representing the marginal likelihood of the observations for a component of the position given a prior mean, in our case the MAP $x_i^*$, is given by
\begin{equation}
\begin{split}
    2 \log  &f(\bm{y}_{x_i}|\bm{t}^y,p) = \\ 
    & -(\bm{y}_{x_i} -x_i^*(\bm{t}^y;p))^T K_{\sigma,i}^{-1} (\bm{y}_{x_i}-x_i^*(\bm{t}^y;p)) \\
    \hspace{12pt} &- \log |K^\sharp_i(\bm{t}^y,\bm{t}^y)| - n_y \log 2 \pi ,
    \label{eq:log_ml}
\end{split}
\end{equation}
\begin{equation}
    \text{where} \quad K_{\sigma,i} :\hspace{-2pt}  = K^\sharp_i(\bm{t}^y,\bm{t}^y) + \sigma^2_{s,i}I.
\end{equation}
Since we model $y_1, \dots, y_{d_z}$ as a priori independent,
\begin{equation}
    f(\bm{y}_x|\bm{t}^y, p)=f(\bm{y}_{x_1}|\bm{t}^y,p) \dots f(\bm{y}_{x_{d_x}}|\bm{t}^y,p).
\end{equation}
The parameter posterior is determined using Bayes rule and Monte Carlo integration of the parameters,
\begin{align}
\begin{split}
    w_p := f(p|\bm{y}_x,\bm{t}^y)
    \approx  \frac{f(\bm{y}_x|\bm{t}^y, p) f(p) (\# \mc{P}_\text{MC})}{\sum_{p^\dagger \in \mc{P}_{\text{MC}}} f(\bm{y}_x|\bm{t}^y,p^\dagger)},
    \label{eq:p_post}
\end{split}
\end{align}
where we assumed $p$ is independent of $\bm{t}^y$.
In the marginalization, $dp$ formally denotes the product measure of the Lebesgue measure for the continuous variables and the counting measure for the discrete variables.
A similar process can be used infer the probability density that an individual parameter $p_i \in \mc{P}_i$ is correct by applying Bayes theorem then marginalizing over all other parameters $p_{\neg i}\in \mc{P}_{\neg i} := \mc{P}_1\times \dots \times \mc{P}_{i-1} \times \mc{P}_{i+1} \times \dots \times \mc{P}_{d_p}$ and again using a priori independence between parameters and with $\bm{t}^y$, 
\begin{align}
\begin{split}
    f(p_i|\bm{y}_x,\bm{t}^y) =& \int_{\mc{P}_{\neg i}} \frac{f(\bm{y}_x|\bm{t}^y, p_{\neg i}, p_i) f(p_i)}{f(\bm{y}_x|\bm{t}^y, p_{\neg i})} f(p_{\neg i})  dp_{\neg i} \\
    \approx& \sum_{p_{\neg i} \in \mc{P}_{\neg i,\text{MC}}} \frac{f(\bm{y}_x|\bm{t}^y, p_{\neg i}, p_i) f(p_i)}{\sum_{p_i^\dagger \in \mc{P}_{i,\text{MC}}} f(\bm{y}_x|\bm{t}^y, p_{\neg i}, p_i^\dagger) },
\label{eq:one_p_post}
\end{split}
\end{align}
where we can borrow $\mc{P}_{\neg i,\text{MC}}$ and $\mc{P}_{i,\text{MC}}$ from the full $\mc{P}_\text{MC}$ due to a priori independence and then the sample count $\# \mc{P}_\text{MC}$ cancels out.
For example, \cref{eq:one_p_post} can be used to find the probability the target will arrive at a given destination or the posterior probability density of turning radii in our Dubins car example.
Finally, we are interested in the posterior probability density of spatial locations given a time and the current observations,
\begin{align}
\begin{split}
    f(x|t,\bm{y}_x,\bm{t}^y) &= \int_{\mc{P}} f(x|t,\bm{y}_x,\bm{t}^y,p) w_p dp \\
    & \approx \sum_{p \in \mc{P}_{\text{MC}}} w_p f(x|t,\bm{y}_x,\bm{t}^y,p)
    \label{eq:x_post}
\end{split}
\end{align}
where the conditional density at a single time is given by
\begin{align}
\begin{split}
    \log f(x_i| & t,\bm{y}_{x_i}, \bm{t}^y,p) = \\
    &-\frac{\log(2 \pi \tilde{\sigma_i}(t;p)^2)}{2} - \frac{(x_i - \tilde{x}_i(t;p))^2}{2 \tilde{\sigma_i}(t;p)^2},
\end{split}
\end{align}which uses the GP posterior mean and variance of the GP corresponding to $p$,
\begin{align}
    \tilde{x}_i(t;p) &= x_i^*(\bm{t}^y;p) + K^\sharp_i(t,\bm{t}^y) K^{-1}_{\sigma,i}(\bm{y}_{x_i} - x_i^*(\bm{t}^y;p)), \nonumber \\
    \tilde{\sigma_i}(t;p)^2 &= k^\sharp_i(t, t) - K^\sharp_i(t,\bm{t}^y)K^{-1}_{\sigma,i} K^\sharp_i(\bm{t}^y,t).
    \label{eq:gp_post}
\end{align}

\section{OPTIMAL RENDEZVOUS PLANNING}
We seek spatiotemporal rendezvous coordinates that maximize the probability of contact under our stochastic model.
A \emph{rendezvous point} $(t, x) \in \mathbb{R}^{d+1}$ is deemed successful when the target lies within distance $R$ at time $t$:
\begin{align}
    \Vert X_{t} - x \Vert \leq R.
\end{align}

Following the model established in the preceding section, we assume the trajectory of the target depends on latent parameters $P$, representing different potential intent, environmental models, or dynamics. Let $\mathcal{P}$ be the set of possible parameters, where $w_p = \mathbb{P}[P = p]$. For each parameter $p \in \mathcal{P}$, let $X_t^{(p)}$ denote the conditional stochastic process. The overall trajectory $X_t$ is then represented by the conditional density:
\begin{equation}
f(x \mid t, p) = f_{X_t \mid P}(x \mid p).
\end{equation}
Our objective is to compute a fixed number of rendezvous points $(t_1, x_1), ..., (t_n, x_n)$ which maximizes the probability of at least one successful rendezvous
\begin{align}
    \max_{(t_1,x_1),...,(t_n, x_n)}\mathbb{P}[\exists i \in \{1,...,n\}\text{: } \Vert X_{t_i} - x_i \Vert \leq R].
\end{align}
By considering the complement (i.e. the probability of all points failing to make contact), this is equivalent to the minimization
\begin{equation}\label{eq:original_min}
\min_{(t_1,x_1), \dots, (t_n, x_n)} \mathbb{P} \left[ \bigcap_{i=1}^n C_i \right],
\end{equation}
where $C_i$ is the event $\Vert X_{t_i} - x_i \Vert > R$. Using the chain rule of probability:
\begin{equation}\label{eq:iterative_min}
\mathbb{P}[C_1, \dots, C_n] = \mathbb{P}[C_1] \prod_{i = 2}^n \mathbb{P}[C_i \mid C_1, \dots, C_{i-1}].
\end{equation}

This decomposition induces a natural sequential decision structure. Computing the globally optimal set of rendezvous points would require joint optimization over all stages, accounting for the dependence of future terms on earlier decisions through posterior updates, which is generally intractable. We therefore adopt a greedy strategy that selects rendezvous points sequentially by minimizing each conditional failure probability in \eqref{eq:iterative_min} under the current belief. After each selection, the distribution of $X_t$ is updated by conditioning on the event that the previous rendezvous attempt failed, yielding a Bayesian update over both parameter weights and trajectory distributions. 
This yields a tractable one-step lookahead policy in belief space: at each stage, the next rendezvous point is chosen by minimizing conditional failure probability under the current posterior, then the belief is updated after an unsuccessful attempt.

\subsubsection{Pursuer Reachability Constraints}

In practical applications, the optimization in \eqref{eq:original_min} cannot be performed over the entire space $\mathbb{R}^{d+1}$. The pursuers are subject to their own kinematic constraints.

We assume the pursuers follow the Dubins car dynamics described in Equation \eqref{eq:car_ode}, though their motion parameters (specifically turning radius and maximum velocity) differ from those of the target. Critically, we consider the regime where the pursuers are significantly slower than the target. In such scenarios, traditional pursuit-evasion strategies or "tail-chasing" maneuvers are fundamentally infeasible, as the pursuer cannot close the distance from behind. In all our numerical experiments, the pursuers are constrained to a constant speed equal to 30\% of the speed of the target, placing the system firmly in a velocity-disadvantaged regime.

Our method overcomes this velocity deficit by leveraging the stochastic model of the target. Instead of reactive pursuit, the problem becomes one of optimal placement: identifying the spatiotemporal "bottlenecks" where the probability of the target's presence is high and the arrival of the pursuer is kinematically guaranteed. This allows for successful rendezvous even when the target possesses a decisive maneuverability advantage.

Let $\mathcal{X} = \{\Gamma^0_1, ..., \Gamma^0_m\}$ be initial configurations representing fixed base-stations or deployment hubs. Each station $\Gamma^0_j$ is capable of launching multiple pursuers, allowing for the simultaneous coverage of several optimal rendezvous points from a single geographic origin. For each $\Gamma^0_j$, let $u_{\Gamma^0_j}(x)$ be the minimum time for the $j$-th pursuer to reach position $x$, calculated as the solution of the Hamilton-Jacobi equation \eqref{eq:HJB_PDE}. Since the value function $u_{\Gamma^0_j}$ is computed once per station, we can evaluate the reachability of any candidate rendezvous point for an entire fleet of pursuers without redundant PDE solves.

The region that is viable for the j-th base-station is then 
\begin{align}
    \mathcal{R}_j = \left\{ (t, x) \in \mathbb{R}^{d+1} \quad \vert \quad u_{\Gamma^0_j}(x) \leq t \right\}.
\end{align}
Beyond simple spatial proximity, the use of the value function $u_{\Gamma^0_j}$ allows for the enforcement of specific contact angles. Since $\Gamma^0_j$ defines the set of initial states (position and orientation), it can naturally encode launch constraints (e.g. runway heading, port exit angle). By incorporating the stochastic heading of the target trajectory $X_t$, we can further refine the reachable set to include only those points where the pursuer arrives at a desired relative orientation. Mathematically, we can describe the constrained reachable set for station $j$ as
\begin{align}
    R_j = \left\{ (t, x) \in \mathbb{R}^{d+1} \hspace{2pt} \big\vert \hspace{2pt} x \in F(t, X_{\cdot}), \hspace{2pt}  u_{\Gamma^0_j}(x) \leq t \right\},
\end{align}
where $F(t, X_{\cdot})$ represents a filter or geometric constraint on the configuration space at time $t$ given the estimated path of the target (typically using the expected position at time $t$). This ensures that the pursuer not only meets the target in space and time but does so in a configuration (e.g., specific approach angle) that maximizes the likelihood of a stable rendezvous.

We define the overall viable set of rendezvous points to be any spatiotemporal point reachable by at least one pursuer
\begin{align}
    \mathcal{R} = \bigcup_{j=1}^m \mathcal{R}_j.
\end{align}

The optimization problem is then constrained to this set, ensuring that for every selected $(t_i, x_i)$, there exists at least one pursuer capable of arriving at the rendezvous point at the same time as the target. Once the optimal points are found, pursuers are assigned to points using a simple bottleneck assignment or by minimizing arrival slack $t_i - u_{\Gamma^0_j}(x_i)$.

\subsubsection{Approximate Belief Update After Rendezvous Failure}

Each unsuccessful rendezvous is incorporated by conditioning the belief on no contact at the selected spatiotemporal point — formally, on $\|X_s - y\| > R$ for attempted rendezvous $(s,y)$. A failure not only informs the distribution at time $s$ but also influences nearby times, since trajectory continuity implies the target was likely absent from the surrounding region as well. Exact posterior recomputation after each such update is intractable in our setting, so we introduce a structured approximate update designed to retain the most important temporal dependencies while remaining computationally efficient.

In our stochastic trajectory model, unsuccessful contacts can arise for two fundamentally different reasons. First, we may have selected an incorrect parameter set, causing the realized trajectory to deviate from the expected one. Alternatively, the parameter set may be correct, and the rendezvous fails purely due to inherent trajectory uncertainty. These two failure modes have markedly different implications: the cross-time correlations inherent to continuous trajectories are informative only when the parameter set is correct, and cannot be leveraged when the failure stems from parameter misidentification. For this reason, we treat the two scenarios separately in our analysis.

First we assume that the parameter set $P$ is correct. In this case we consider the implication of the spatiotemporal correlation of the continuous trajectories to update our distribution estimates $f(x|t, p)$ for the trajectories of each parameter set $p$. For a failed rendezvous point $(s, y)$, the updated trajectory distribution is computed via
\begin{align}
    f(x \mid \Vert X_s - y \Vert > R, t, p) =  \frac{f(x|t, p) g(x, y)}{\int f(z|t, p) g(z, y) dz}, \\
    g(z, z') = 1 - \exp\!\left(-\frac{\|z - z'\|^2}{2\sigma(R)^2}\right).
\end{align}

This governs sequential belief updates for trajectories under parameter $p$. The parameter weights are updated as follows, again assuming $(s, y)$ is a failed rendezvous attempt:
\begin{align}
    \mathbb{P}[P=a \vert \Vert X_s - y \Vert > R] = \frac{\mathbb{P}[\Vert X^{(a)}_s - y \Vert > R]}{\mathbb{P}[\Vert X_s - y \Vert > R]} \mathbb{P}[P = a].
\end{align}

Together, these updates iteratively refine the trajectory distribution after each failed rendezvous, consistent with the sequential structure in \eqref{eq:iterative_min}.
In practice, the pursuer’s dynamics impose constraints on which regions of the spatiotemporal domain can be reached from the initial state. Depending on the starting position, velocity, and admissible controls, large portions of the domain may be physically unattainable. To account for this, we pre-compute the reachable set $\mathcal{R} \subseteq \mathbb{R}^{d+1}$, defined as the collection of spatiotemporal points the pursuers can reach. All probability computations and optimization steps are then restricted to this subdomain. Points outside $\mathcal{R}$ are masked, ensuring that the algorithm considers only physically realizable rendezvous opportunities.

The proposed strategy is summarized in Algorithm~\ref{alg:optimization}. By leveraging the Bayes updates derived above, the planner iteratively selects the most promising reachable spatiotemporal points, while accounting for the evolving belief over the target's trajectory and motion parameters.

%

\begin{algorithm}[h]
\caption{MAP-Based Reachable Rendezvous Planner}
\label{alg:optimization}
\begin{algorithmic}[1] 
    \Require Prior mixture weights $\{w_p\}$, conditional densities $\{f(t,x|p)\}$, reachable set $\mathcal{R}$, and number of rendezvous attempts $n$.
    \State Discretize domain into grid $\mathcal{G}$
    \State Mask points $(t, x) \notin \mathcal{R}$
    \For{$k=1, \dots, n$} \Comment{Select $n$ rendezvous points}
        \State \textbf{Step 1: Selection}
        \State Compute current failure probability:
        \State \quad $P(t, x) = \sum_{p \in \mathcal{P}_{MC}} w_p \int_{\|x-y\| > R} f(t, y|p) dy$
        \State Select point $(s_k, y_k) \in \mathcal{R} \cap \mathcal{G}$ that minimizes $P(t, x)$

        \State \textbf{Step 2: Belief Update (Conditioning on Failure)}
        \State Update mixture weights of parameter uncertainty:
        \State \quad$w^\prime_p \gets w_p \int_{\|x - y\| > R} f(x|s, p) dx$ 
        \State \quad $w_p \gets w_p^\prime / \sum_{p \in \mathcal{P}_{MC}} w_p^\prime$ 
        \State Update conditional densities of trajectories:
        \State \quad $f(x|t, p) \gets \frac{f(x|t, p) g(x, y)}{\int f(z|t, p) g(z, y) dz}$ 
        \State \textbf{Record} rendezvous point $(s, y)$ 
    \EndFor
\end{algorithmic}
\end{algorithm}

This procedure efficiently identifies a small set of high-value rendezvous points while accounting for uncertainty in both target dynamics and parameter configurations.

\section{NUMERICAL RESULTS}

\begin{figure*}[t]
\centering
\begin{subfigure}{0.8\linewidth}
    \includegraphics[width=\textwidth]
    {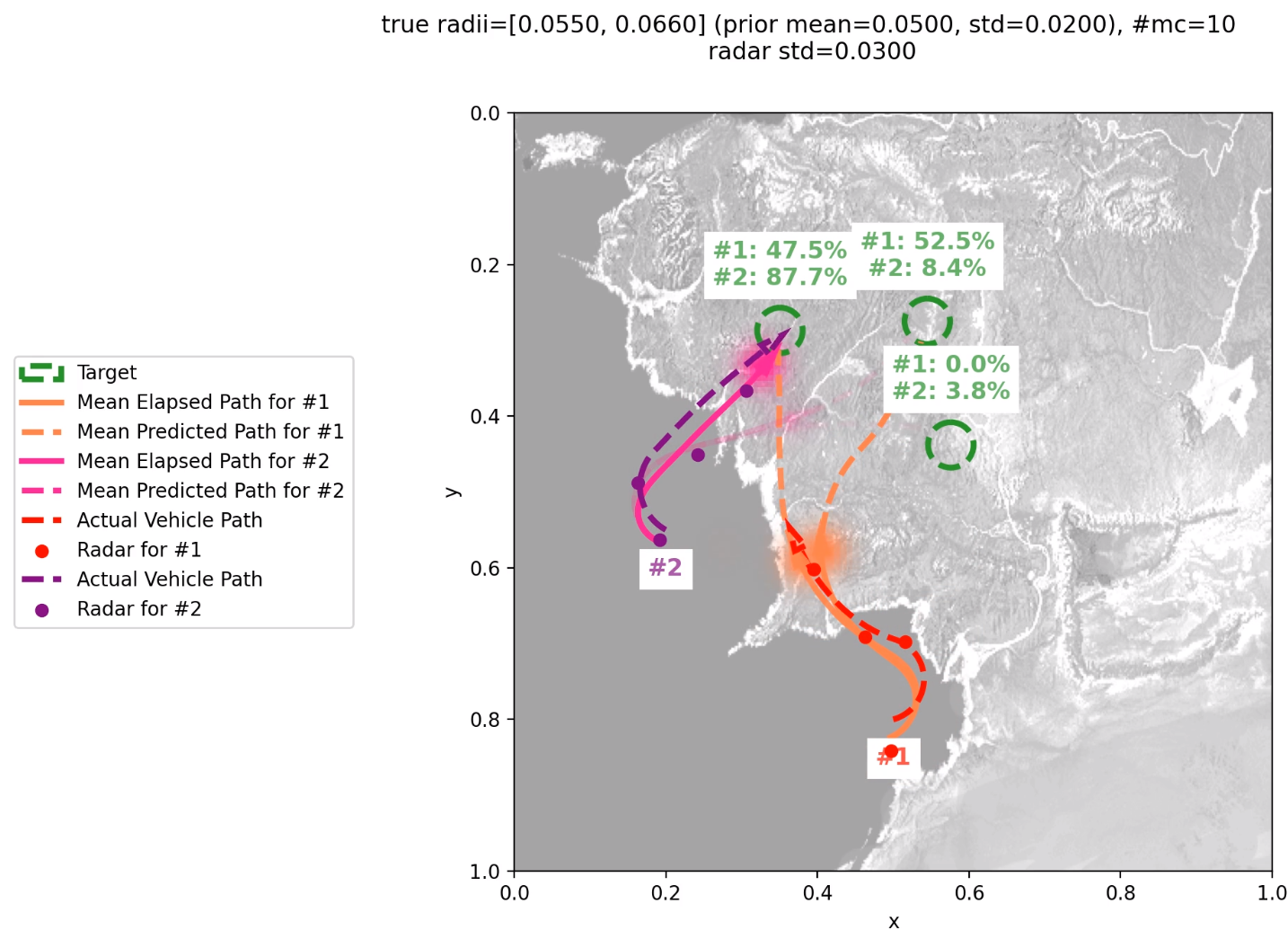}
\end{subfigure}
\caption{Multiple target, multiple destination estimates given very noisy data}
\label{fig:gp_paths}
\end{figure*}

We present a scenario using our Dubins car example (\cref{fig:gp_paths}), where target entity \#1 has true turning radius 0.055 and \#2 has 0.066. The density of $x$ is computed from \cref{eq:x_post} for each target, and the hit probability from \cref{eq:one_p_post}. Observations $y$ of position in the unit domain and orientation (in radians) are i.i.d. Gaussian with standard deviation $\sigma_s = 0.03$. Both targets move at known unit speed toward the leftmost destination region, where all regions have radius 0.03. The turning radius and destination index jointly define the true parameter $p$.

Bayesian estimation uses a Gaussian prior with mean 0.05 and standard deviation 0.02 on the turning radius $\rho$, and a uniform prior over the destination index $i$, forming $f(p)$. For the solver, we set $\beta = \sigma_s$, corresponding to the equivalent weight in standard GPR. For MAP estimation, we use kernel $k(t,t')=\exp(-\tfrac{1}{4}0.05^{-2}(t-t')^2)$, while GPR uses $k^\sharp(t,t')=\sigma_s^2 \exp(-\tfrac{1}{4}0.1^{-2}(t-t')^2)$.

\begin{figure}
\centering
\begin{subfigure}{1.0\linewidth}
    \includegraphics[width=\textwidth]
    {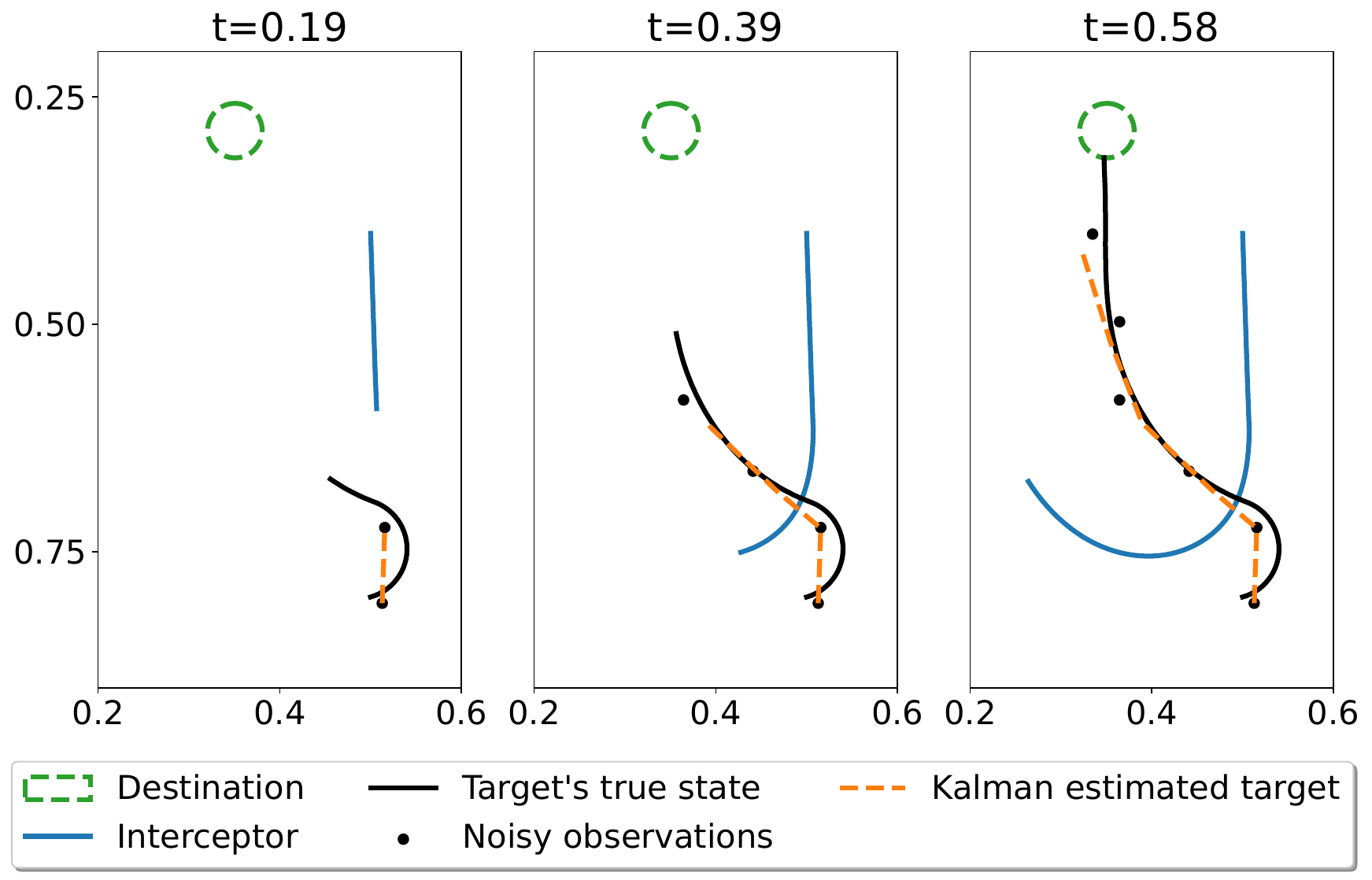}
\end{subfigure}
\caption{Given sparse observations, proportional guidance using Kalman filtering misses on evenly matched target and pursuer entities due to lag caused by lack of a nonlinear dynamics constraint in the filter}
\label{fig:follower}
\end{figure}

\begin{figure}
\centering
\begin{subfigure}{1.0\linewidth}
    \includegraphics[width=\linewidth]
    {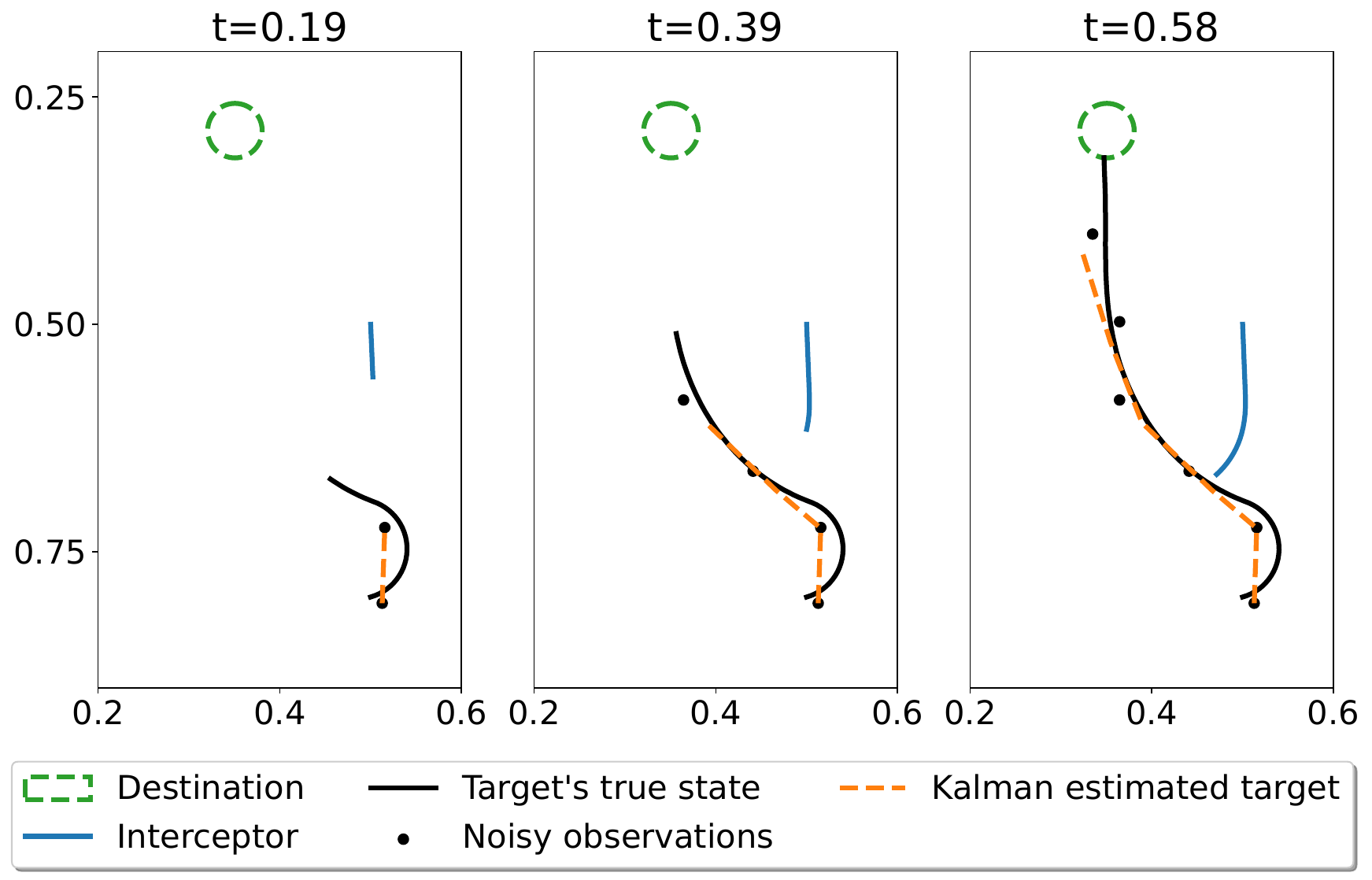}
\end{subfigure}
\caption{In the same slow-pursuer scenario, follower dynamics using a Kalman filter utterly fails due to lack of planning}
\label{fig:slow_compare}
\end{figure}

Despite sparse observations, the Bayesian inference correctly identifies the destination and turning radius. As shown in \cref{fig:follower}, this outperforms existing methods, where Kalman filtering and proportional guidance fail to yield a kinematically viable rendezvous when the pursuer is subject to significant velocity disadvantages. The proportional controller uses the same minimum turning radius as the target in \cref{fig:follower}, while \cref{fig:slow_compare} shows a case with a pursuer at 30\% of the target’s speed.

We assume that the pursuers follow the dynamics stated in equation \eqref{eq:car_ode} with a turning radius of $\rho = 0.05$. Further we restrict the rendezvous points to make contact at a perpendicular angle (i.e., from the side). This perpendicular-rendezvous constraint significantly reduces the feasible rendezvous set, as the pursuer must arrive not only at the correct location and time, but also with a precise relative orientation. We consider two different scenarios for pursuer deployment:
\begin{itemize}
    \item Unrestricted Scenario: The initial configuration of the pursuers is entirely unconstrained.
    \item Restricted Scenario: A more constrained setup where the initial configuration must have a specified angle (limited to heading east or west).
\end{itemize}

\begin{figure}[h!]
    \centering
    \includegraphics[width=0.49\linewidth]{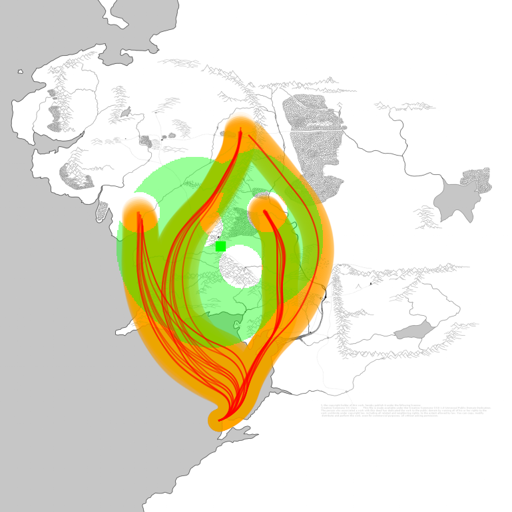}
    \includegraphics[width=0.49\linewidth]{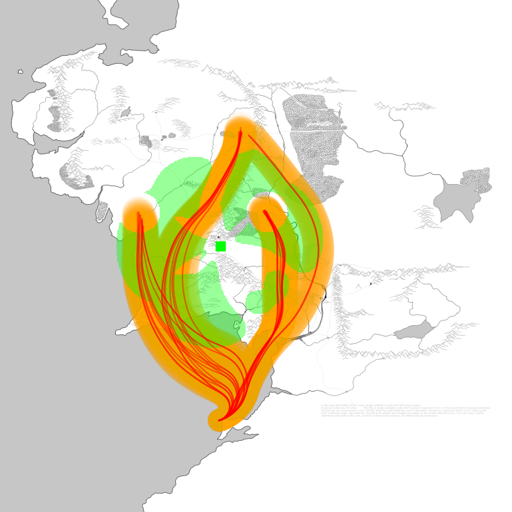}
    \caption{Reachable regions $\mathcal{R}$ at terminal time for the two pursuer deployment scenarios targeting perpendicular contact. Left: The unrestricted scenario where the initial pursuer configuration is unconstrained. Right: The restricted scenario where initial headings are limited to east or west. The green dot denotes the launch base station, and the light green regions highlight spatiotemporal coordinates where perpendicular contact is feasible according to the value function $u$.}
    \label{fig:reachable_region_terminal}
\end{figure}

These two scenarios yield differing reachable regions $\mathcal{R}$, with the restricted scenario having fewer feasible choices than the unrestricted case. The set $\mathcal{R}$ at terminal time for both scenarios is illustrated in figure~\ref{fig:reachable_region_terminal}. The green dot indicates the base station from which pursuers are launched and the light green region shows which locations can be reached at the correct time and angle.

While $\mathcal{R}$ generally expands as the time horizon increases, granting pursuers more time to navigate to distant locations, the reachable region is not necessarily strictly growing. Because of the specific maneuvers required to achieve a perpendicular rendezvous, certain spatiotemporal gaps may exist where a location is temporarily unreachable before becoming feasible again at a later time.

To compute the reachable regions, we associate every location with the closest mean trajectory point $\mathbb{E}[X_t^{(p)}]$, compute the corresponding perpendicular rendezvous angle, and check reachability at that angle using the value function $u$ at every time point.

\begin{figure}[h!]
    \centering
    \includegraphics[width=0.49\linewidth]{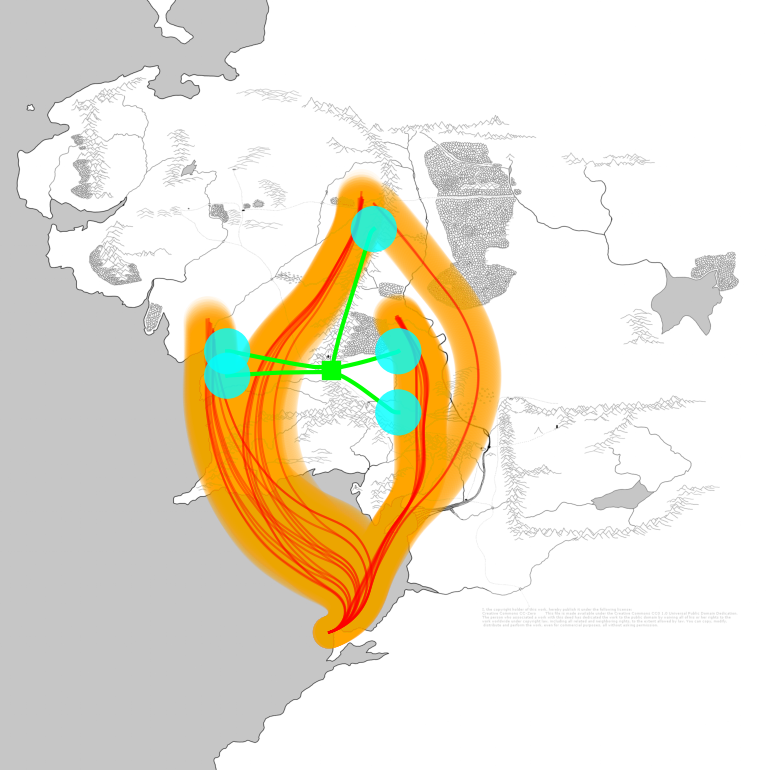}
    \includegraphics[width=0.49\linewidth]{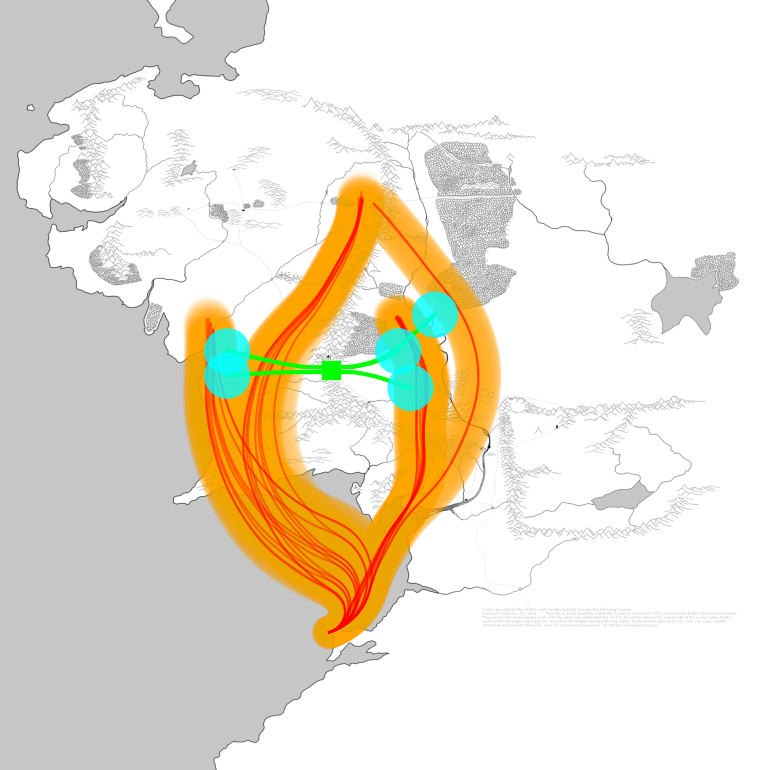}
    \caption{Optimal rendezvous points and pursuer trajectories for unrestricted (left) and restricted (right) scenarios. The green lines represent the optimal rendezvous paths derived from the value function $u$, targeting the stochastic trajectory ensembles (orange tubes) and their mean paths (red). The terminal points are marked by cyan circles of radius $R$, representing the optimized spatiotemporal coordinates for successful rendezvous and parameter disambiguation.}
    \label{fig:intercept_paths_terminal}
\end{figure}

The outcome of our optimization procedure, as governed by Algorithm~\ref{alg:optimization}, is shown in Figure~\ref{fig:intercept_paths_terminal}. To determine the optimal rendezvous coordinates, we analyze the stochastic trajectory ensembles (represented by orange tubes), where each ensemble corresponds to a predicted path distribution for a specific parameter set. The red curves denote the mean trajectories under each model, and the surrounding orange regions capture the associated spatial uncertainty.

Algorithm~\ref{alg:optimization} selects the spatiotemporal positions where a rendezvous attempt is most informative according to our sequential objective. These targeted coordinates are marked by the cyan circles of radius $R$. As shown, these points are positioned to maximize the likelihood of intersecting at least one plausible trajectory while simultaneously providing strong disambiguation between parameter sets.

The computed rendezvous points provide us with the specific time and spatial coordinates to target. Using the value function $u$, which is the solution of equation~\eqref{eq:HJB_PDE} and was previously computed for the construction of the reachable region $\mathcal{R}$, we can compute the optimal rendezvous paths. This is achieved by applying the optimal control law:
\begin{align}
    \frac{d\theta}{dt}(x)= -\frac{\nabla_\theta u(x)}{\Vert \nabla_\theta u(x) \Vert}.
\end{align}
We note that by selecting the initial location of the pursuer and computing optimal trajectories in reverse from the computed rendezvous point to its initial position, the underlying Hamilton–Jacobi equation only needs to be solved once per pursuer initial location.
The resulting optimal paths are illustrated in Figure~\ref{fig:intercept_paths_terminal} as green lines terminating at the cyan target regions.

In practice, pursuers are not necessarily launched simultaneously, but are deployed only when required by the optimal solution. Because reachability is encoded directly in the value function $u$, the resulting trajectories implicitly determine both the launch time and the path of each pursuer. The green trajectories correspond to pursuers en route to their designated rendezvous points, following their respective optimal controls toward future spatiotemporal targets.

\section{CONCLUSION}
In this paper, we presented a robust probabilistic framework for planning rendezvous points for a dynamically constrained, slower pursuer and a target observed through sparse and noisy data. By combining a kernel-based MAP estimation with a Gaussian process correction, the proposed method effectively models the spatial uncertainty of unknown target trajectories. Sequential optimization via greedy minimization of failure probability over the current belief state provides a systematic approach to predictive planning under uncertainty.

In the scenarios considered here, our results indicate advantages over the Kalman-filter-plus-proportional-guidance baseline. Simulations with Dubins car dynamics show that Kalman filtering combined with proportional guidance fails to produce adequate pursuer trajectories, due to lag and a lack of predictive planning, particularly when the pursuer has nonlinear dynamics or lower speed. In contrast, our Bayesian inference framework rapidly and accurately infers critical target parameters, such as the destination region and turning radius, despite the sparsity of the observational data.

Furthermore, we demonstrated the framework's adaptability to strict kinematic and operational constraints. By leveraging the Hamilton-Jacobi-Bellman value function, the method accurately computes reachable regions and identifies feasible spatiotemporal rendezvous points, even under severe restrictions such as mandatory perpendicular rendezvous angles or constrained initial launch configurations. In our examples, the algorithm converts stochastic trajectory estimates into feasible rendezvous plans and associated pursuer trajectories. It identifies spatiotemporal coordinates with high predicted rendezvous value under the current belief and yields corresponding launch times and steering controls through the precomputed value function.

Ultimately, this sequential optimization and Bayesian updating strategy overcomes the inherent limitations associated with slower pursuers and strict kinematic constraints.
Future work includes optimizing the kernel hyperparameters for both the MAP and GPR steps, which is challenging for real-time scenarios; evaluating our method on more complex systems of ODEs; and extending our planning algorithm to non-myopic objectives.
By providing a reliable method for generating optimal, executable rendezvous trajectories under uncertain dynamics, this framework offers broad utility for real-world autonomous rendezvous applications, including unmanned aerial vehicle refueling, spacecraft servicing, autonomous surface vessel operations, search and rescue missions, and missile defense.

\section{ACKNOWLEDGMENTS}

A portion of this research was funded by Lockheed Martin.

\FloatBarrier 

\bibliographystyle{unsrt} 
\bibliography{references}

\end{document}